\documentclass[english,12pt,oneside]{amsproc}
\usepackage[T2A]{fontenc}
\usepackage[english]{babel}
\usepackage{sseq}
\usepackage{graphics}
\usepackage{amsfonts, amssymb, amscd, amsmath}
\usepackage{latexsym}
\usepackage[matrix,arrow,curve]{xy}
\usepackage{mathabx,mathtools}
\usepackage{color}
\usepackage{mathrsfs}
\usepackage{mathdots}
\usepackage{pigpen}
\usepackage{tikz}
\usetikzlibrary{matrix}

 \DeclareMathOperator{\cone}{Cone}

\DeclareMathOperator{\ver}{Vert}

\DeclareMathOperator{\soc}{Soc}

\newcommand{\Zo}{\mathbb{Z}}
\newcommand{\Ro}{\mathbb{R}}
\newcommand{\Co}{\mathbb{C}}
\newcommand{\Qo}{\mathbb{Q}}

\newcommand{\ko}{\Bbbk}

\newcommand{\simc}{\!\!\sim}

\newcommand{\br}{\widetilde{\beta}}

\newcommand{\pd}{\texttt{PD}}
\newcommand{\Ins}{I_{N\!S}}

\newcommand{\eqd}{\stackrel{\text{\tiny def}}{=}}

\newcommand{\minel}{\hat{0}}

\newcommand{\wh}[1]{{\widehat{#1}}}

\newcommand{\E}[1]{(E_{#1})}
\newcommand{\dif}[1]{(d_{#1})}

\newcommand{\Hr}{\widetilde{H}}
\newcommand{\dd}{\partial}

\newcommand{\Qh}{\wh{Q}}
\newcommand{\Xh}{\wh{X}}

\newcounter{stmcounter}[section]
\newcounter{thcounter}

\numberwithin{equation}{section}

\theoremstyle{plain}

\newtheorem{thm}[thcounter]{Theorem}

\newtheorem{defin}[stmcounter]{Definition}

\theoremstyle{definition}

\newtheorem{rem}[stmcounter]{Remark}
\newtheorem{con}[stmcounter]{Construction}

\begin{document}

\title{Topological model for h''-vectors of simplicial manifolds}

\author[Anton Ayzenberg]{Anton Ayzenberg}
\address{Department of Mathematics, Osaka City University, Sumiyoshi-ku, Osaka 558-8585, Japan.}
\email{ayzenberga@gmail.com}

\date{\today}
\thanks{The author is supported by the JSPS postdoctoral fellowship program.}

\begin{abstract}
Any manifold with boundary gives rise to a Poincare duality
algebra in a natural way. Given a simplicial poset $S$ whose
geometric realization is a closed orientable homology manifold,
and a characteristic function, we construct a manifold with
boundary such that graded components of its Poincare duality
algebra have dimensions $h_k''(S)$. This gives a clear topological
evidence for two well-known facts about simplicial manifolds: the
nonnegativity of $h''$-numbers (Novik--Swartz theorem) and the
symmetry $h''_k=h''_{n-k}$ (generalized Dehn--Sommerville
relations).
\end{abstract}

\maketitle


%
%
%
%
%
%
%

\section{Introduction}\label{SecIntro}

Let $S$ be a pure simplicial poset of dimension $n-1$ and
$[m]=\{1,\ldots,m\}=\ver(S)$ be the set of its vertices. Let $\ko$
be a ground ring which is either a field or~$\Zo$. A map
$\lambda\colon [m]\to \ko^n$ is called a \emph{(homological)
characteristic function} if, for any maximal simplex $I\in S$, the
set of vertices of $I$ maps to the basis of $\ko^n$. We suppose
that there is a fixed basis in $\ko^n$, and, for any vertex $i\in
[m]$, the value $\lambda(i)$ has coordinates
$(\lambda_{i,1},\ldots,\lambda_{i,n})$ in this basis.

Let $\ko[S]$ be the face ring of $S$ (see \cite{StPos,BPposets}).
By definition, $\ko[S]$ is a commutative associative graded
algebra over $\ko$ generated by formal variables $v_I$, one for
each simplex $I\in S$, with relations
\[
v_{I_1}\cdot v_{I_2}=v_{I_1\cap I_2}\cdot\sum_{J\in I_1\vee
I_2}v_J,\qquad v_{\emptyset}=1.
\]
Here $I_1\vee I_2$ denotes the set of least upper bounds of
$I_1,I_2\in S$, and $I_1\cap I_2\in S$ is the intersection of
simplices (it is well-defined and unique when $I_1\vee
I_2\neq\emptyset$). The summation over an empty set is assumed to
be $0$. For topological reasons we take the doubled grading on the
ring: the generator $v_I$ has degree $2|I|$, where $|I|$ denotes
the rank of $I$. The natural map $\ko[m]=\ko[v_1,\ldots,v_m]\to
\ko[S]$ defines the structure of $\ko[m]$-module on $\ko[S]$.

Any characteristic function $\lambda\colon[m]\to \ko^n$ determines
the set of linear elements:
\[
\theta_1=\sum_{i\in [m]}\lambda_{i,1}v_i,\quad \theta_2=\sum_{i\in
[m]}\lambda_{i,2}v_i,\quad\cdots,\quad\theta_n=\sum_{i\in
[m]}\lambda_{i,n}v_i\in \ko[S]
\]
(these elements have degree $2$, but we will use the term
``linear'' when its meaning is clear from the context). The
definition of characteristic function implies that
$\theta_1,\ldots,\theta_n$ is a \emph{linear system of parameters}
in $\ko[S]$ (e.g.\cite[Lm.3.5.8]{BPnew}). Moreover, any linear
system of parameters arise from some characteristic function in
this way. Let $\Theta$ be the ideal in $\ko[S]$ generated by the
elements $\theta_1,\ldots,\theta_n$.

The quotient $\ko[S]/\Theta$ is a finite-dimensional vector space.
The standard reasoning in commutative algebra implies that,
whenever $S$ is Cohen--Macaulay, the dimension of the homogeneous
component $(\ko[S]/\Theta)_{2k}$ is $h_k$, the $h$-number of $S$
\cite{StPos}.

When $S$ is Buchsbaum, the additive structure of $\ko[S]/\Theta$
is still independent of the choice of characteristic function but
dimensions of homogeneous components have more complicated
description. By Schenzel's theorem \cite{Sch,NS} the dimension of
$(\ko[S]/\Theta)_{2k}$ is
\[
h_k'\eqd h_k+{n\choose
k}\left(\sum_{j=1}^{k-1}(-1)^{k-j-1}\br_{j-1}(S)\right),
\]
where $\br_k(S)=\dim \Hr_k(S;\ko)$.

Recall that the socle of a $\ko[m]$-module $\mathcal{M}$ is a
$\ko$-subspace
\[
\soc\mathcal{M} \eqd \{y \in \mathcal{M} \mid \ko[m]^{+}\cdot y =
0\},
\]
where $\ko[m]^+$ is the maximal graded ideal of the ring $\ko[m]$.
Since the products with polynomials of positive degrees are
trivial, the socle is a $\ko[m]$-submodule of $\mathcal{M}$.
In~\cite{NS} Novik and Swartz proved the existence of certain
submodules in $\soc(\ko[S]/\Theta)$ for any Buchsbaum simplicial
poset. Namely, in degree $2k<2n$ there exists a vector subspace
\[
(\Ins)_{2k}\subseteq\soc(\ko[S]/\Theta)_{2k},
\]
isomorphic to ${n\choose k}\Hr^{k-1}(S;\ko)$ (this notation means
the direct sum of ${n\choose k}$ copies of $\Hr^{k-1}(S;\ko)$).
Let $\Ins$ denote the direct sum of $(\Ins)_{2k}$ over all $k$,
where we assumed $(\Ins)_{2n}=0$. Since $\Ins$ lies in the socle,
it is a $\ko[m]$-submodule. Moreover, $\Ins$ is an ideal in
$\ko[S]/\Theta$ (for simplicial complex this fact easily follows
from the surjectivity of the map $\ko[m]\to\ko[S]$, and for
simplicial poset, whose geometrical realization is a homology
manifold, this was checked in \cite[Rem. 8.3]{AyV3}). Thus we may
consider the quotient ring $(\ko[S]/\Theta)/\Ins$. The dimension
of its homogeneous component of degree $2k$ is equal to $h_k''$
where
\[
h_k'' \eqd h_k'-{n\choose k}\br_{k-1}(S) = h_k+{n\choose
k}\left(\sum_{j=1}^{k}(-1)^{k-j-1}\br_{j-1}(S)\right),
\]
for $0\leqslant k\leqslant n-1$, and $h_n''=h_n'$. In particular,
$h_k''\geqslant 0$ for any Buchsbaum simplicial poset.

Now we restrict to the case when the ground ring is either $\Zo$
or $\Qo$. The class of Cohen--Macaulay simplicial posets contains
an important subclass of sphere triangulations. By abuse of
terminology we call simplicial poset a homology sphere (resp.
manifold) if its geometrical realization is a homology sphere
(resp. manifold).

Every homology sphere is Cohen--Macaulay. For homology spheres the
ring $\ko[S]/\Theta$ is a Poincare duality algebra (this is not
surprising in view of Danilov--Jurkiewicz and Davis--Januszkiewicz
theorems). In general one can prove this by the following
topological argument. Consider the cone over $|S|$ endowed with a
dual simple face stratification and consider the identification
space $X_S=(\cone|S|\times T^n)/\simc$, similar to the
construction of quasitoric manifold \cite{DJ}. Using the same
ideas as in \cite{DJ}, one can prove that the cohomology algebra
of $X_S$ over $\ko$ is isomorphic to $\ko[S]/\Theta$ (see e.g.
\cite{MasPan}). When $\ko=\Zo$, the space $X_S$ is a homology
manifold over integers. In case $\ko=\Qo$, this space is a
homology manifold over $\Qo$. In both cases the Poincare duality
over the corresponding ring implies that $\ko[S]/\Theta$ is a
Poincare duality algebra. In particular, this proves
Dehn--Sommerville relations for homology spheres: $h_k=h_{n-k}$.

The goal of this paper is to construct a topological model for
homology manifolds. Every closed homology manifold is a Buchsbaum
simplicial poset, so the ideal $\Ins\subset \ko[S]/\Theta$ is
defined. There holds

\begin{thm}[cf.{\cite[Th.1.4]{NSgor}}]\label{thmPDA}
Let $S$ be a simplicial poset whose geometrical realization is a
closed connected orientable homology manifold. When the ground
ring is either $\Qo$ or $\Zo$, the ring $(\ko[S]/\Theta)/\Ins$ is
a Poincare duality algebra.
\end{thm}

It was proved in \cite{NSgor} that $(\ko[S]/\Theta)/\Ins$ is a
Gorenstein ring, which implies the theorem over any field when $S$
is a simplicial complex. The theorem gives a straightforward
evidence for the generalized Dehn--Sommerville relations for
manifolds: $h''_{k}=h''_{n-k}$.

The idea of our proof is the following. In Section
\ref{SecPoinPrelim} we associate a Poincare duality algebra with
any manifold with boundary $(M,\dd M)$, either smooth,
topological, or homological. This algebra will be denoted
$\pd^*_{(M,\dd M)}$. Given any homology manifold $S$, instead of
taking the cone (as in the case of spheres) we consider the
topological space $\Qh=|S|\times [0,1]$. This space is a manifold
with two boundary components: $\dd_0\Qh$ and $\dd_1\Qh$. Consider
the identification space $\Xh=(\Qh\times T^n)/\simc$ with the
relation collapsing certain torus subgroups over the points of
$\dd_0\Qh$ similar to a quasitoric case, and not touching the
points over $\dd_1\Qh$. The space $\Xh$ is a homology manifold
with boundary; its boundary consists of points over $\dd_1\Qh$.
Then Theorem \ref{thmPDA} is an immediate consequence~of

\begin{thm}\label{thmCohomOfCollarModel}
The algebra $\pd^*_{(\Xh,\dd\Xh)}$ is isomorphic to
$(\ko[S]/\Theta)/\Ins$.
\end{thm}

The only place in the arguments, where we need the restriction on
a ground ring, is the construction of torus spaces. The relation
$\sim$ collapses certain compact subgroups of the compact torus
$T^n$, and this identification cannot be defined for general
characteristic functions over general fields. Nevertheless, if the
characteristic function $\lambda$ over $\ko$ can be represented as
$\lambda'\otimes\ko$ for some characteristic function $\lambda'$
over $\Zo$ (or $\Qo$), then the statements hold true over a field
$\ko$ and this particular choice of characteristic function.

%
%
%
%
%
%
%

\section{Poincare duality algebras}\label{SecPoinPrelim}

\begin{defin}
A finite-dimensional, graded, associative, graded-commutative,
connected algebra $A^*=\bigoplus_{k=0}^dA^k$ over $\ko$ is called
Poincare duality algebra of formal dimension $d$, if
\begin{enumerate}
\item $A^d\cong \ko$;

\item The product map $A^k\otimes A^{d-k}\stackrel{\times}{\rightarrow} A^d$
is a non-degenerate pairing for all $k=0,\ldots,d$. Over integers
the finite torsion should be mod out.
\end{enumerate}
\end{defin}

While the motivating examples of Poincare duality algebras are
cohomology of connected orientable closed manifolds, there exist
another natural source of duality algebras.

\begin{con}\label{conPDmdfbndr}
Let $(M,\dd M)$ be a compact connected orientable manifold with
boundary, $\dim M=d$. As a technical requirement we will also
assume that $M$ contains a neighborhood of $\dd M$ of the form
$\dd M\times[0,\varepsilon]$. Consider the $\ko$-module
$A^*=\bigoplus_{k=0}^dA^k$, where
\[
A^k=\begin{cases} H^0(M), \mbox{ if } k=0;\\
\mbox{image of } \imath^*\colon H^k(M,\dd M)\to H^k(M),\mbox{ if
}0<k<d;\\
H^d(M,\dd M),\mbox{ if }k=d.
\end{cases}
\]
The homomorphism $\imath^*\colon H^k(M,\dd M)\to H^k(M)$ is
induced by the inclusion $\imath\colon(M,\emptyset)\hookrightarrow
(M,\dd M)$.

There is a well-defined product on $A^*$ induced by the
cup-products in cohomology. Indeed, let $a_1\in A^{k_1}$ and
$a_2\in A^{k_2}$. If either $k_1$ or $k_2$ is zero, then there is
nothing to define, since $A^0$ is spanned by the unit of the ring.
If $k_1+k_2<d$ then $a_1\cdot a_2$ is just the product of two
elements in the ring $H^*(M)$. This product lies in the image of
$H^*(M,\dd M)$. If $k_1+k_2=d$, then we may consider the elements
$b_1,b_2\in H^*(M,\dd M)$ such that
$\imath^*(b_\epsilon)=a_{\epsilon}$, and take their product in the
ring $H^*(M,\dd M)$. This gives an element in $H^d(M,\dd M)= A^d$
which we call the product of $a_1$ and $a_2$. It is easily seen
that this element does not depend on the choice of representatives
$b_1,b_2$ for the elements $a_1,a_2$.

The Poincare--Lefschetz duality implies that the pairing between
$A^k$ and $A^{d-k}$ is non-degenerate \cite{Br}. Thus $A^*$ is a
Poincare duality algebra. We denote it by $\pd_{(M,\dd M)}$ and
call \emph{the Poincare duality algebra of a manifold with
boundary}.
\end{con}

\begin{rem}\label{remHomolPD}
By Poincare--Lefchetz duality, instead of cohomology we can work
with homology. We have
\[
\pd_{(M,\dd M)}^k\cong\begin{cases} H_d(M,\dd M), \mbox{ if } k=0;\\
\mbox{image of } \imath_*\colon H_{d-k}(M)\to H_{d-k}(M,\dd
M),\mbox{ if
}0<k<d;\\
H_0(M),\mbox{ if }k=d,
\end{cases}
\]
and the product is given by the intersection product in homology.
\end{rem}

%
%
%
%
%
%
%

\section{Collar model}\label{SecCollar}

\subsection{Buchsbaum simplicial posets}
Recall that a finite partially ordered set (poset) $S$ is called
\emph{simplicial} if: (1)~there exists a unique minimal element
$\minel\in S$; (2)~for each element $J\in S$, the lower order
ideal $\{I\in S\mid I\leqslant J\}$ is isomorphic to the poset of
faces of a $k$-simplex for some number $k$. Elements of $S$ are
called simplices, and the atoms of $S$ are called vertices. The
number $k$ is called the dimension of $I$ and the number $|I|=\dim
I+1$, which is equal to the number of vertices of $I$, is called
the rank of $I$. In the following we assume that $S$ is pure of
dimension $n-1$, which means that all maximal simplices of $S$
have the same rank $n$. Let $S'$ be the barycentric subdivision of
$S$. For each proper simplex $I\in S\setminus\minel$ consider the
following subsets of the geometrical realization $|S|=|S'|$:
\[
G_I=|\{(I_0< I_1<\ldots)\in S'\mbox{ such that } I_0\geqslant
I\}|,
\]
\[
\dd G_I=|\{(I_0< I_1<\ldots)\in S'\mbox{ such that } I_0>I\}|.
\]
The subset $G_I$ is called the face of $|S|$ dual to $I$. A
simplicial poset $S$ is called \emph{Buchsbaum} (over $\ko$) if
$H_j(G_I,\dd G_I;\ko)=0$ for any $I\in S\setminus\minel$ and
$j\neq\dim G_I$. In particular, any homology manifold is
Buchsbaum, since in this case $G_I$ are homological cells.

\subsection{Collar model}
Consider the compact $n$-torus with a fixed coordinate
representation $T^n=\{(t_1,\ldots,t_n)\mid t_s\in \Co,|t_s|=1\}$.
If $\ko$ is either $\Zo$ or $\Qo$, the vector
$w=(w_1,\ldots,w_n)\in \ko^n$ determines a compact 1-dimensional
subgroup
$t^w=\{(e^{2\pi\sqrt{-1}w_1t},\ldots,e^{2\pi\sqrt{-1}w_nt})\mid
t\in \Ro\}$. Let $[m]=\ver(S)$ be the set of vertices of $S$ and
let $\lambda\colon[m]\to \ko^n$ be a characteristic function over
$\Zo$ or $\Qo$. Let $T_i\subset T^n$ denote the one dimensional
subgroup $t^{\lambda(i)}$. For a simplex $I\in S$, let $T_I$
denote the product of the one-dimensional subgroups $T_i$
corresponding to the vertices of $I$ (the product is taken inside
$T^n$). The definition of characteristic function implies that
$T_I$ is a compact subtorus of $T^n$ of dimension $|I|$.

Consider the space $\Qh=|S|\times [0,1]$ which will be called the
\emph{collar} of $|S|$. Let $\dd_\epsilon\Qh$ denote the subset
$|S|\times\{\epsilon\}$ for $\epsilon=0,1$. The faces $G_I$ can be
considered as the subsets of $\dd_0\Qh\subset \Qh$. To make the
notation uniform, we set $G_{\minel}=\Qh$ and
$T_{\minel}=\{1\}\subset T^n$.

\begin{con} Consider the identification space $\Xh=(\Qh\times
T^n)/\simc$, where the points $(x,t)$, $(x',t')$ are identified
whenever $x=x'\in G_I$ and $t^{-1}t'\in T_I$ for some simplex
$I\in S$. Let $f\colon \Qh\times T^n\to \Xh$ denote the quotient
map, and $\mu$ denote the projection to the first factor,
$\mu\colon \Xh\to\Qh$. The preimage $\mu^{-1}(G_I)$ is denoted by
$X_I$. Let $\dd_1\Xh$ denote the subset $\dd_1\Qh\times T^n\subset
\Xh$.
\end{con}

\subsection{Absolute and relative spectral sequences}
The dual face structure on $|S|$ induces a topological filtration
\[
Q_0\subset Q_1\subset \ldots\subset Q_{n-1}=\dd_0\Qh\subset Q_n=
\Qh,
\]
which lifts to the orbit type filtration on $\Xh$:
\[
X_0\subset X_1\subset \ldots\subset X_{n-1}\subset X_n=\Xh.
\]
Let $\E{\Qh}^1_{p,q}=H_{p+q}(Q_p,Q_{p-1})\Rightarrow H_{p+q}(\Qh)$
and $\E{\Xh}^1_{p,q}=H_{p+q}(X_p,X_{p-1})\Rightarrow H_{p+q}(\Xh)$
be the homological spectral sequences associated with these
filtrations. By the result of \cite{AyV1}, whenever $S$ is
Buchsbaum, the map
\[
f_*^2\colon \bigoplus_{q_1+q_2=q}\E{\Qh}^2_{p,q_1}\otimes
H_{q_2}(T^n)\to \E{\Xh}^2_{p,q}
\]
is an isomorphism for $p>q$ and injective for $p=q$. The rightmost
column of $\E{\Qh}^*_{*,*}$ (which is the source of all higher
differentials) vanish: $\E{\Qh}^1_{n,*}\cong
H_{n+*}(\Qh,\dd_0\Qh)=0$, since the collar $\Qh$ collapses to
$\dd_0\Qh$. Similar for $\E{\Xh}^*_{*,*}$. Thus there are no
higher differentials $d^{\geqslant 2}$ in both spectral sequences.

We also need the homological spectral sequences for the relative
homology:
\[
\E{(\Qh,\dd_1\Qh)}^r_{p,q}\Rightarrow
H_{p+q}(\Qh,\dd_1\Qh)=0,\qquad
\E{(\Xh,\dd_1\Xh)}^r_{p,q}\Rightarrow H_{p+q}(\Xh,\dd_1\Xh).
\]
The first pages are the following:
\[
\E{(\Qh,\dd_1\Qh)}^1_{p,q}=\begin{cases}H_{p+q}(Q_p,Q_{p-1}),\mbox{
if }p<n;\\
H_{n+q}(\Qh,Q_{n-1}\sqcup\dd_1\Qh),\mbox{ if } p=n,
\end{cases}
\]
\[
\E{(\Xh,\dd_1\Xh)}^1_{p,q}=\begin{cases}H_{p+q}(X_p,X_{p-1}),\mbox{
if }p<n;\\
H_{n+q}(\Xh,X_{n-1}\sqcup\dd_1\Xh),\mbox{ if } p=n.
\end{cases}
\]
Note that the rightmost terms $\E{(\Qh,\dd_1\Qh)}^1_{n,q}$ have
the form:
\[
H_{n+q}(\Qh,\dd_0\Qh\sqcup\dd_1\Qh)\cong H_{n+q}(|S|\times[0,1],
|S|\times\{0,1\})\cong H_{n+q-1}(S),
\]
and the higher differentials
\[
\dif{(\Qh,\dd_1\Qh)}^r\colon \E{(\Qh,\dd_1\Qh)}^*_{n,-r+1}\to
\E{(\Qh,\dd_1\Qh)}^*_{n-r,0}\cong H_{n-r}(S)
\]
are isomorphisms (so that the spectral sequence for
$(\Qh,\dd_1\Qh)$ collapses to zero). Similar to the non-relative
case, the induced map
\[
f_*^2\colon
\bigoplus_{q_1+q_2=q}\E{(\Qh,\dd_1\Qh)}^2_{p,q_1}\otimes
H_{q_2}(T^n)\to \E{(\Xh,\dd_1\Xh)}^2_{p,q}
\]
is an isomorphism for $p>q$ and injective for $p=q$ (this follows
from the general method developed in \cite{AyV1}).

\subsection{Proof of Theorem \ref{thmCohomOfCollarModel}}

The proof essentially relies on calculations made in \cite{AyV3}.
If $S$ is a connected orientable homology manifold, then $\Xh$ is
a connected orientable homology manifold with the boundary
$\dd_1\Xh\cong |S|\times T^n$. The boundary admits a collar
neighborhood as required in Construction \ref{conPDmdfbndr}. For
$I\neq\minel$ the subset $X_I$ is a closed submanifold of
codimension $2|I|$ lying in the interior of $\Xh$. It is called
the \emph{face submanifold}, and its homology class $[X_I]\in
H_{2n-2|I|}(\Xh)$ is called the \emph{face class}. Note that for
$|I|\neq 0$ the classes $[X_I]$ appear in the spectral sequence
$\E{\Xh}^*_{*,*}$ as the free generators of the group:
$\E{\Xh}^1_{q,q}$ with $q=n-|I|$. The relations on these classes
in $H_{2q}(\Xh)$ are precisely the images of first differentials,
hitting the group $\E{\Xh}^1_{q,q}$. In \cite[Prop.4.3 and
Lm.8.2]{AyV3} we checked that these relations are the same as the
linear relations on $v_I$ in the quotient ring
$(\ko[S]/\Theta)_{2(n-q)}$ when $q\leqslant n-2$. If $q=n-1$,
there are no relations on $[X_I]\in H_{2n-2}(\Xh)$ since there are
no differentials hitting the group $\E{\Xh}^1_{n-1,n-1}$.

In addition to face classes, there exist other homology classes in
$H_*(\Xh)$, namely the classes coming from the part of the
spectral sequence below the diagonal. They lie in the groups
$\E{\Xh}^2_{p,q}\cong H_p(\Qh)\otimes H_q(T^n)$ for $q<p<n$. In
\cite{AyV3} we called them \emph{spine classes}.

Let us keep track on the behavior of homology classes, when they
map to the relative homology by the homomorphism $\imath_*\colon
H_*(\Xh)\to H_*(\Xh,\dd_1\Xh)$. Again, we may look at their
representatives in the spectral sequence $\E{(\Xh,\dd_1\Xh)}^*$.
At this time, higher differentials are nontrivial. All spine
classes of $H_*(\Xh)$ are killed by higher differentials. Indeed,
they lie in the part of the relative spectral sequence which is
isomorphic to $\E{(\Qh,\dd_1\Qh)}^*_{*,*}\otimes H_*(T^n)$, and
the latter sequence collapses to $0$.

On the other hand, the diagonal cells of the relative spectral
sequence are hit by higher differentials as well. Thus there are
more relations on $[X_I]$ in the group $H_*(\Xh,\dd_1\Xh)$ than in
the group $H_*(\Xh)$. The higher differential
\[
\dif{(\Xh,\dd_1\Xh)}^r\colon \E{(\Xh,\dd_1\Xh)}^1_{n,n-2r+1}\to
\E{(\Xh,\dd_1\Xh)}^1_{n-r,n-r}
\]
is injective, and gives an inclusion of $H_{n-r}(S)\otimes
H_{n-r}(T^n)$ into $\E{(\Xh,\dd_1\Xh)}^1_{n-r,n-r}$. Under the
degree reversing identification $[X_I]\leftrightarrow v_I$ (and by
Poincare duality in~$S$), this inclusion gives the Novik--Swartz
submodule $(\Ins)_{2r}\cong {n\choose r}\Hr^{r-1}(S)$ inside
$(\ko[S]/\Theta)_{2r}$, for $r\geqslant 2$ (see details in
\cite{AyV3}). When $r=1$, only the first differential
$\dif{(\Xh,\dd_1\Xh)}^1$ hits the cell
$\E{(\Xh,\dd_1\Xh)}^1_{n-1,n-1}$. Its image corresponds to the
linear span of $\theta_1,\ldots,\theta_n$ in $\ko[S]_{2}$. There
is no Novik--Swartz submodule $\Ins$ in $(\ko[S]/\Theta)_{2}$
since $S$ is connected.

These considerations prove that the image of the map
\[
\imath_*\colon H_{2n-2r}(\Xh)\to H_{2n-2r}(\Xh,\dd_1\Xh)
\]
is isomorphic to the homogeneous component of
$(\ko[S]/\Theta)/\Ins$ of degree $2r$ for each $0<r<n$.

When $r=n$, the submodule $(\Ins)_{2n}$ is trivial. Thus
$H_0(\Xh)$ coincides with
$(\ko[S]/\Theta)_{2n}=((\ko[S]/\Theta)/\Ins)_{2n}$. The group
$H_{2n}(\Xh,\dd_1\Xh)\cong \ko$ is obviously identified with
$((\ko[S]/\Theta)/\Ins)_0\cong\ko$.

Theorem \ref{thmCohomOfCollarModel} now follows from Remark
\ref{remHomolPD} and the fact that the correspondence
$[X_I]\leftrightarrow v_I$ translates the intersection product on
$\Xh$ to the product in the face ring.

\end{document}